\documentclass[12pt]{article}

\usepackage{graphicx}
\usepackage{amsmath}
\usepackage{amsthm}
\usepackage{amssymb}

\usepackage{natbib}

\theoremstyle{definition}
\newtheorem{thm}{Theorem}
\newtheorem{cor}{Corollary}
\newtheorem{lem}{Lemma}
\newtheorem{rem}{Remark}

\theoremstyle{definition}
\newtheorem{defn}{Definition}
\newtheorem{exmp}{Example}

\newcommand{\mbR}{{\mathbb{R}}}
\newcommand{\mbZ}{{\mathbb{Z}}}

\begin{document}

\title{
A Jacobian inequality for gradient maps on the sphere
and its application to directional statistics
}
\author{Tomonari SEI}

\maketitle

\begin{abstract}
In the field of optimal transport theory,
an optimal map is known to be a gradient map of
a potential function satisfying cost-convexity.
In this paper, the Jacobian determinant
of a gradient map is shown to be log-concave
with respect to a convex combination of the potential functions
when the underlying manifold is the sphere and the cost function
is the distance squared.
The proof uses the non-negative cross-curvature property
of the sphere recently established by Kim and McCann, and Figalli and Rifford.
As an application to statistics,
a new family of probability densities on the sphere
is defined in terms of cost-convex functions.
The log-concave property of the likelihood function
follows from the inequality.
\end{abstract}

\setlength{\baselineskip}{18pt}

\section{Introduction}

In recent years, the theory of optimal transport
has been actively studied.
In particular, properties of the optimal transport map
on Riemannian manifolds are well established.
The existence and uniqueness theorem for the optimal transport map
on Riemannian manifolds
was proved by \cite{mccann2001}; this result extended the pioneering work of \cite{brenier1991}
for the Euclidean case.
He showed that optimal transport is given by the gradient map
of a so-called cost-convex function.
On the other hand, for statistical data analysis on Euclidean space, it
is useful to consider convex combinations of convex functions
in order to construct various probability density functions
(\cite{sei2006}, \cite{sei2009b}).
In this paper, we show that when the underlying space is the sphere,
the convex combination of cost-convex functions is actually cost-convex
(Lemma~\ref{lem:convex}) and the Jacobian determinant of
the resultant gradient map is log-concave with respect to the
convex combination (Theorem~\ref{thm:Jacobian}).
This result is an extension of the Jacobian interpolation inequality
shown by \cite{cordero2001}.
We refer to our Jacobian inequality as \textit{the Jacobian inequality} throughout this  paper,
for simplicity.

Our result is related to the regularity theory of optimal transport maps.
Here we consider some recent studies in this field.
\cite{ma2005} showed that regularity of the transport map for general cost functions
on Euclidean space is assured if a geometrical quantity
called the \textit{cost-sectional curvature} is positive.
Conversely, \cite{loeper2005} showed 
that non-negativity of the cost-sectional curvature
is necessary for regularity.
He also showed that non-negativity of the cost-sectional curvature implies
non-negativity of the usual sectional curvature
if the cost function is the squared distance on a Riemannian manifold.
However, the converse does not hold (\cite{kim2007counter}).
Comprehensive assessment on the theory of optimal transport
has been published (\cite{villani2009}).
A relevant concept is \textit{the cross-curvature} (\cite{kim2007}).
\cite{kim2008} and \cite{FigalliRifford2009} 
independently showed that the sphere $S^n$ has almost positive cross-curvature.
In general, the cost-sectional curvature is non-negative
if the cross-curvature is non-negative.
In the present paper, we use the non-negative cross-curvature property of the sphere
to prove our main results.

We show that our Jacobian inequality
opens several doors for applications to directional statistics.
In this field,
a family of probability densities is used
to analyze given directional data, such as locations on the earth.
For example, a test on the directional character of given data
is constructed via families of probability density functions on the sphere.
Directional statistics has a long history
since \cite{fisher1953}
and a comprehensive text on this subject
has been published (\cite{mardia2000}).

We define a probability density function on the sphere
by the gradient maps of cost-convex functions.
Although, in the context of optimal transport,
one usually considers push-forward of probability densities,
we construct a family of densities by means of \textit{pull-back}
of probability densities.
This follows from the fact that a pull-back density has an explicit expression
for the likelihood function needed for statistical analysis.
The density function does not need any special functions
such as the modified Bessel function, which usually appear in directional statistics.
Furthermore, the Jacobian inequality
implies that the likelihood function is log-concave with respect to the
statistical parameters.
This property is reasonable for computation of the maximum likelihood estimator.
We propose more specific models
and show graphical images of each probability density.
In terms of analysis of real data,
we present the result of density estimation for some astronomical data.

This paper is organized as follows.
In Section~\ref{sec:main}, we present basic notation
and state our main theorem.
In Section~\ref{sec:appl}, we construct a family of
probability density functions on the sphere
and apply them to directional statistics.
All mathematical proofs of the main theorem and lemmas
are given in Section~\ref{sec:proof}.
Finally we present a discussion in Section~\ref{sec:discussion}.

\section{Main theorem} \label{sec:main}

Let $S^n$ be the $n$-dimensional unit sphere.
The tangent space at $x\in S^n$ is denoted by $T_xS^n$.
The geodesic distance (arc length) between $x$ and $y$ in $S^n$
is denoted by $d(x,y)$.
The cost function is $c(x,y)=(1/2)d(x,y)^2$.
If one uses Euclidean coordinates in $\mbR^{n+1}$
to express $S^n$,
then $d(x,y)=\cos^{-1}(x^{\top}y)$,
where the range of $\cos^{-1}$ is $[0,\pi]$.

The \textit{$c$-transform} $\phi^c$ of a function $\phi:S^n\to\mbR$
is defined by
\begin{align}
 \phi^c(y)\ =\ \sup_{x\in S^n}
 \left\{-c(x,y)-\phi(x)\right\}.
 \label{eqn:c-transform}
\end{align}
The function $\phi$ is said to be cost-convex, or \textit{$c$-convex}, if
$(\phi^c)^c=\phi$.
Examples of $c$-convex functions will be given
in Section~\ref{sec:appl}.
By compactness of $S^n$,
a function $\phi$ is $c$-convex if and only if
for any $x\in S^n$
there exists some (not necessarily unique) $y\in S^n$
such that $c(x,y)+\phi(x)=\inf_{z\in S^n}\{c(z,y)+\phi(z)\}$.

The image of the exponential map of $v\in T_xS^n$ at $x\in S^n$,
denoted by $\exp_x(v)$,
is the end point of the geodesic starting
at $x$ with the initial vector $v$.
More explicitly,
if one uses Euclidean coordinates in $\mbR^{n+1}$
to express $S^n$ and $T_xS^n$,
the exponential map is written as
$\exp_x(v)=(\cos|v|)x+(\sin|v|)(v/|v|)$,
where $|v|$ denotes the Euclidean norm of the vector $v$.
The exponential map $\exp_x$ is a diffeomorphism from $\{v\in T_xS^n\mid |v|<\pi\}$
to $S^n\setminus\{x'\}$, where $x'$ is the antipodal point of $x$.

The following lemma is a consequence of
the non-negative cross-curvature property of the sphere
established by \cite{kim2008} and \cite{FigalliRifford2009}.
See Section~\ref{sec:proof} for a proof.

\begin{lem}[Convex combination of $c$-convex functions] \label{lem:convex}
 If $\phi_0$ and $\phi_1$ are $c$-convex,
 then for each $t\in[0,1]$ the function $(1-t)\phi_0(x)+t\phi_1(x)$ of $x$
 is also $c$-convex.
\end{lem}

\begin{rem}
 \cite{FKM2009} showed Lemma~\ref{lem:convex} simultaneously and independently from us.
 Indeed, they showed more general result, in that
 the convexity of the space of $c$-convex functions
 is necessary and sufficient condition for the
 non-negative cross-curvature property (see Theorem 3.2 in \cite{FKM2009}).
\end{rem}

We define $G_{\phi}(x)=\exp_x(\nabla\phi(x))$
as long as $\phi$ is differentiable at $x$,
where $\nabla$ is the gradient operator.
Following \cite{delanoe2006},
we call $G_{\phi}:S^n\to S^n$ the \textit{ gradient map}
associated with the \textit{ potential function} $\phi$.
The map $G_{\phi}$ is differentiable at $x$
if $|\nabla\phi(x)|<\pi$ and $\phi$ has its Hessian at $x$.
It is known that any $c$-convex $\phi$ on any compact Riemannian manifold
is Lipschitz and therefore differentiable almost everywhere.
Furthermore, $\phi$ has a Hessian almost everywhere in the Alexandrov sense,
and therefore $G_{\phi}(x)$ is differentiable almost everywhere
(see \cite{mccann2001} and \cite{cordero2001}).
These technical facts on differentiability
are important for the theory of optimal transport.
However, we will not need them because, for statistical applications,
we can assume from the beginning that $G_{\phi}(x)$
is differentiable except at a finite set of
points (see Section~\ref{sec:appl}).

For any $c$-convex functions $\phi_0$ and $\phi_1$,
by Lemma~\ref{lem:convex},
the convex combination $\phi_t(x)=(1-t)\phi_0(x)+t\phi_1(x)$
is $c$-convex.
We define an interpolation of gradient maps by
\[
 F_t(x)
 \ =\ G_{\phi_t}(x)
 \ =\ \exp_x(\nabla\phi_t(x)),
 \quad t\in[0,1].
\]
Assume that for each $i\in\{0,1\}$, $|\nabla\phi_i(x)|<\pi$ and $\phi_i(x)$ has its Hessian at $x$.
Then it is easy to see that, for any $t\in[0,1]$, $|\nabla\phi_t(x)|<\pi$ and $\phi_t(x)$ has 
its Hessian defined at $x$.
We define the Jacobian determinant $J_t(x)=\mathrm{Jac}(F_t(x))=\det(dF_t/dx)$
with respect to any orthonormal basis
on $T_xS^n$ and $T_{F_t(x)}S^n$ with suitable orientations.

The following theorem is our main result.

\begin{thm}[Jacobian inequality] \label{thm:Jacobian}
 Let $\phi_0$ and $\phi_1$ be two $c$-convex functions.
 Let $x$ be a point in $S^n$ such that, for each $i=0,1$, $|\nabla\phi_i(x)|<\pi$
 and $\phi_i$ has its Hessian defined at $x$.
 Then the Jacobian determinant $J_t(x)$ defined above
 is log-concave with respect to $t$.
 It is equivalent to the inequality
 \[
 \log J_t(x)
 \ \geq\ 
 (1-t)\log J_0(x)
 +t\log J_1(x),\quad t\in[0,1].
 \]
 We refer to the above inequality as the \textit{Jacobian inequality}
 in this paper.
\end{thm}

\begin{rem}
 This theorem is an extension
 of the result obtained by \cite{cordero2001}.
 They showed a similar inequality
 under the additional assumption that $\phi_0\equiv 0$,
 as a corollary of a stronger inequality
 related to the geometric-arithmetic inequality.
 It is not known whether the stronger one
 holds for our case $\phi_0\not\equiv 0$
 (see also Remark~\ref{rem:Cordero}).
 \qed
\end{rem}

\section{Application to directional statistics} \label{sec:appl}

\subsection{Probability densities induced by gradient maps}

In \cite{sei2006} and \cite{sei2009b}, the author proposed a family of
probability density functions in terms of gradient maps on
Euclidean space, where a probability density is constructed
as a pull-back of some fixed measure (typically Gaussian) pulled by a 
gradient map.
The notion can be directly extended to probability density functions
on the sphere.

For statistical application,
we will consider only $c$-convex functions $\phi$ such that
the gradient map $G_{\phi}$ is an isomorphism on $S^n$
and $\phi$ has its Hessian defined everywhere except at a finite set of points.
We define some related terminology.

\begin{defn}[Wrapping potential function] \label{defn:wrapping}
 We say that a function $\phi$ is a \textit{ wrapping potential function}
 if $\phi$ is $c$-convex, $\phi$ has its Hessian defined everywhere except for a finite set of points
 and $G_{\phi}$ is an isomorphism on $S^n$.
 Let $W(S^n)$ be the set of all wrapping potential functions.
\end{defn}

We have the following lemma.

\begin{lem} \label{lem:convex-Gamma}
 If $\phi_0$ and $\phi_1$ are in $W(S^n)$,
 then the interpolation $\phi_t=(1-t)\phi_0+t\phi_1$ ($t\in[0,1]$)
 is also in $W(S^n)$.
\end{lem}

We construct a probability density function for each $\phi\in W(S^n)$.
Let $U$ be a random variable on $S^n$ distributed uniformly.
Then, since $x\mapsto G_{\phi}(x)$ is bijective,
we can define a random variable on $S^n$
by $X=G_{\phi}^{-1}(U)$.
The probability density function of $X$
with respect to the uniform measure
is $p_{\phi}(x)=\mathrm{Jac}(G_{\phi}(x))$,
where the symbol Jac refers to the Jacobian determinant.
In other words, we define $p_{\phi}(x)$
by the \textit{pull-back} measure of the uniform measure
pulled by the gradient map $G_{\phi}$.

At this point, we describe the exact sampling method
of the probability density function $p_{\phi}(x)$.
A sampling procedure is important if one needs
to calculate expectations by the Monte Carlo method.
From the definition, it is clear that the random variable
$X=G_{\phi}^{-1}(U)$ with a uniformly random variable $U$ on $S^n$
has density $p_{\phi}(x)$.
Hence if we can generate $U$ and solve the equation $G_{\phi}(X)=U$ effectively,
we obtain a random sample $X$.
Indeed, $U$ is quite easily generated, for example,
by normalization of a standard Gaussian sample in $\mbR^{n+1}$.
To solve $G_{\phi}(X)=U$,
it is sufficient to find
the unique minimizer of the
function $c(x,U)+\phi(x)$ with respect to $x$
since the following lemma holds.

\begin{lem} \label{lem:exact}[Lemma 7 of \cite{mccann2001}]
 If $\phi$ is $c$-convex and $u=G_{\phi}(x_0)$
 is defined at $x_0\in S^n$, then the
 unique minimizer of $c(x,u)+\phi(x)$ with respect to $x$ is $x_0$.
\end{lem}

\begin{figure}[htbp]
 \begin{center}
  \begin{tabular}{cc}
   \includegraphics[width=6.5cm,clip]{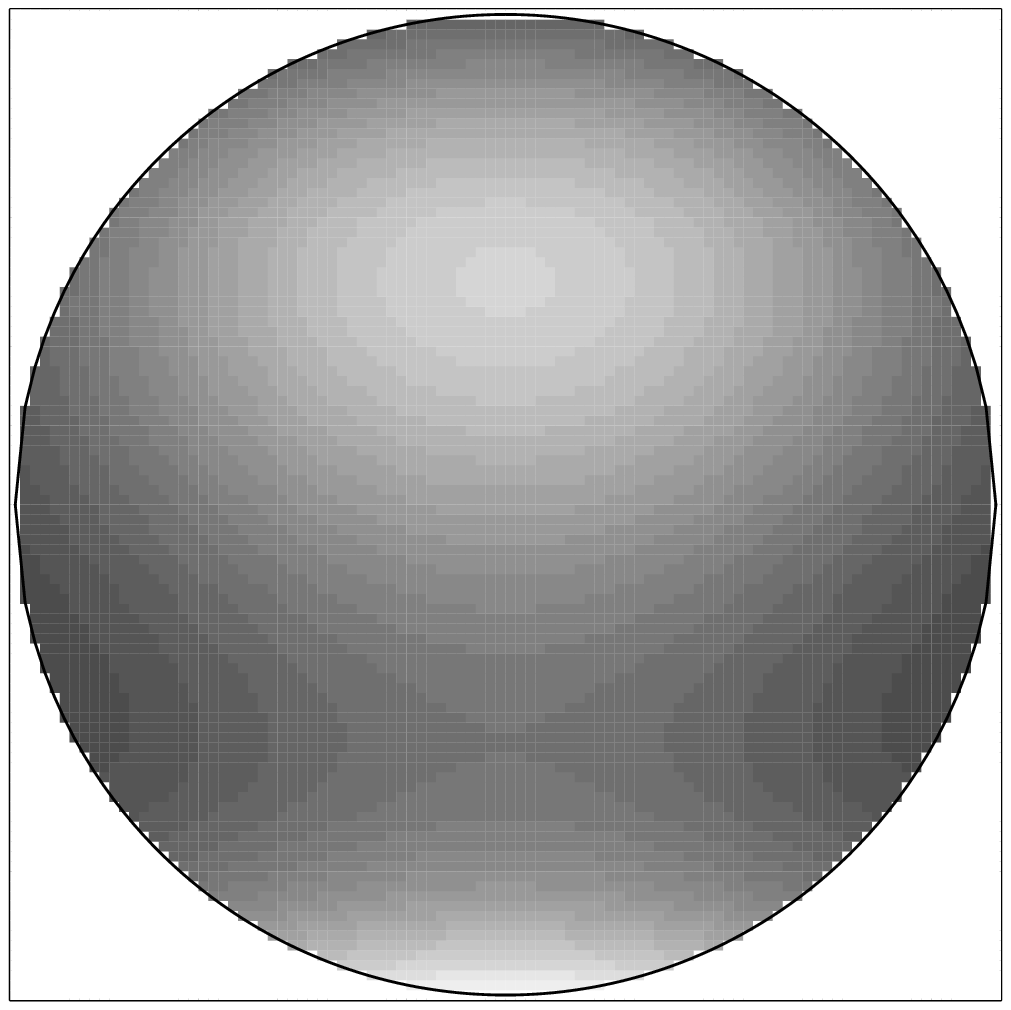} &
   \includegraphics[width=6.5cm,clip]{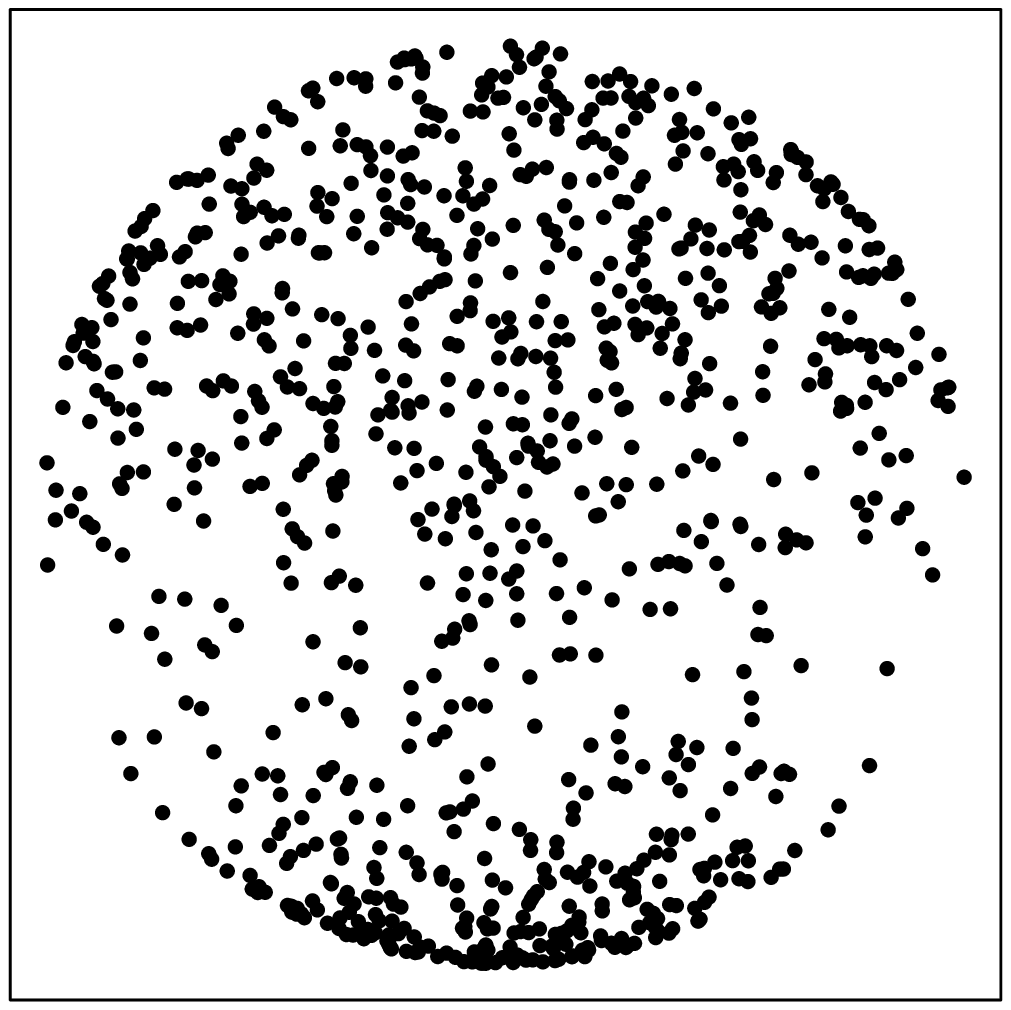}\\
   (a) A density function on the sphere. &
   (b) Samples. \\
  \end{tabular}
 \end{center}
 \caption{
 Exact sampling.
 (a) a density function (a white region indicates high density)
 and (b) 2000 sampled data points.
 The $c$-convex function used is
 $\phi(x)=0.5\cos(2d(x,e_1))+0.5\cos(3d(x,e_2))$ for $x\in S^2$,
 where $e_1$ and $e_2$ denote unit vectors
 along the horizontal and vertical axes, respectively.
 See Subsection~\ref{subsec:examples} for details.
 Only the northern hemisphere is drawn.
 The number of points on the northern hemisphere was 967 in this experiment.
 The program code was written in R and the computational time
 for sampling was about ten seconds.
 }
 \label{fig:sample}
\end{figure}

Thus our task is to solve the (deterministic) minimization problem.
Although the minimization problem of $c(x,U)+\phi(x)$
is not convex in the usual sense,
the objective function has no local minimum, by $c$-convexity.
Hence the problem is efficiently solved by generic optimization packages.
An example of sampling is illustrated in Figure~\ref{fig:sample}.

\subsection{Spherical gradient model}

We consider a finite-dimensional set of probability densities on the sphere.
In statistics, a finite-dimensional set of probability densities
is called a \textit{ statistical model}.
An unknown parameter $\theta$ that parameterizes the density functions
is estimated from observed data points $x(1),\ldots,x(N)\in S^n$.
One of the most important estimators is the maximum likelihood estimator
that maximizes the likelihood function $\prod_{t=1}^Np(x(t)|\theta)$
with respect to $\theta$.

We construct a new statistical model using $c$-convex functions.
Recall that the set $W(S^n)$ of wrapping potential functions
is a convex space (Lemma~\ref{lem:convex-Gamma}).
We can consider a finite-dimensional subspace as follows.
Let $\phi_{(i)}\in W(S^n)$ for $i=1,\ldots,p$.
Define
\begin{align*}
 \phi_{\theta}(x)
 \ =\ \sum_{i=1}^p\theta_i\phi_{(i)}(x),
\end{align*}
where $\theta=(\theta_i)_{i=1}^p$ ranges over
a convex subset $\Theta$ of $\mbR^p$ such that
$\phi_{\theta}\in W(S^n)$ for any $\theta\in\Theta$.
By Lemma~\ref{lem:convex-Gamma} and the elementary fact that $0\in W(S^n)$,
we can use the simplex $\{\theta\mid \theta_i\geq 0,\sum_{i=1}^p\theta_i\leq 1\}$
as $\Theta$.
Let $p(x|\theta)$ be the probability density function
induced by $\phi_{\theta}(x)\in W(S^n)$,
that is,
\begin{align}
 p(x|\theta)
 \ =\ \mathrm{Jac}(G_{\phi_{\theta}}(x)).
 \label{eqn:gradient-model}
\end{align}
We call the family (\ref{eqn:gradient-model})
\textit{the spherical gradient model}.

The maximum likelihood estimator for the spherical gradient model (\ref{eqn:gradient-model})
is reasonably computed
by the following corollary of Theorem~\ref{thm:Jacobian}.

\begin{cor} \label{cor:statistics}
 Define $p(x|\theta)$ by (\ref{eqn:gradient-model}).
 Then, for any data points $x(1),\ldots,x(N)\in S^n$,
 the likelihood function $\prod_{k=1}^Np(x(k)|\theta)$
 is log-concave with respect to $\theta$.
\end{cor}

\begin{rem}
 As an anonymous referee pointed out,
 in Euclidean space, there are results on
 \textit{convexity along generalized geodesics}
 (Chapter 9 of \cite{ambrosio2005}; see also \cite{villani2009}).
 Here a generalized geodesic is defined by
 the set of measures \textit{pushed forward} by
 the gradient maps $\{G_{\phi_t}\}_{t\in[0,1]}$.
 On the other hand, we consider \textit{pull-back measures}
 in this paper as, for example, (\ref{eqn:gradient-model}) indicates.
\end{rem}

\subsection{Examples} \label{subsec:examples}

We give some examples of the spherical gradient model (\ref{eqn:gradient-model}).
Recall $d(x,y)$ denotes the length between $x$ and $y$ on $S^n$.
All the examples are combinations of \textit{rotationally symmetric} functions
$f(d(x,z))$, where $z\in S^n$ and
$f\in C^2([0,\pi])$.
The $k$-th derivative of $f$ is denoted by $f^{(k)}$.
The following lemma is fundamental.

\begin{lem} \label{lem:rotationally-symmetric}
 Assume that $f^{(1)}(0)=f^{(1)}(\pi)=0$ and $f^{(2)}(r)>-1$
 for almost all $r\in[0,\pi]$.
 Then for each $z\in S^n$ the function $f(d(x,z))$ of $x$ is
 in $W(S^n)$.
\end{lem}

Let $\mathcal{F}$ be the set of functions on $[0,\pi]$
that satisfy the assumption in Lemma~\ref{lem:rotationally-symmetric}.
Choose $p$ pairs $\{(f_i,z_i)\}_{i=1}^p$ from $\mathcal{F}\times S^n$.
Then we can define the spherical gradient model (\ref{eqn:gradient-model}) with
\begin{align}
 \phi_{\theta}(x)
 \ =\ \sum_{i=1}^p
 \theta_if_i(d(x,z_i))
 \quad \theta=(\theta_i)_{i=1}^p\in\Theta,
 \label{eqn:rot-sym-potential}
\end{align}
where $\Theta$ is a convex subset of $\mbR^p$
such that $\phi_{\theta}\in W(S^n)$ for all $\theta\in\Theta$.

\begin{rem} \label{rem:rotationally}
 If $p=1$, the resultant density
 $p(x|\theta)$ is a function of $d(x,z)$
 for some $z\in S^n$.
 In directional statistics, such a probability density function
 is called rotationally symmetric.
 \qed
\end{rem}

We briefly touch on known distributions on the sphere in statistics.
A very well-known distribution on the sphere is the
\textit{ von Mises-Fisher distribution} defined by
\begin{align}
  p(x|\mu)
  \ =\ \left(\frac{|\mu|}{2}\right)^{(n+1)/2}
  \frac{1}{\Gamma((n+1)/2)I_{(n+1)/2-1}(|\mu|)}\exp(\mu^{\top}x)
  \label{eqn:vonMises-Fisher}
\end{align}
in Euclidean coordinates of $\mbR^{n+1}$,
where $\mu\in\mbR^{n+1}$
and $I_{\nu}$ denotes the modified Bessel function
of the first kind and order $\nu$.
A more general distribution is the
\textit{ Fisher-Bingham distribution}
defined by
\begin{align}
  p(x|\mu,A)
  \ =\ \frac{1}{a(\mu,A)}
  \exp\left(\mu^{\top}x+x^{\top}Ax\right),
  \label{eqn:Fisher-Bingham}
\end{align}
where $a(\mu,A)$ is a normalizing factor to ensure that
$\int p(x|\mu,A)dx=1$.
See \cite{mardia2000} for details.

We return to our spherical gradient model (\ref{eqn:gradient-model})
with (\ref{eqn:rot-sym-potential}).
The following explicit formula due to a general
expression (\ref{eqn:Jacobian-determinant}) is useful for practical implementation:
\begin{align*}
 p(x|\theta)
 &\ =\ 
 \left(\sin|v_{\theta}|/|v_{\theta}|\right)^{n-1}
 \det\left(
 xx^{\top}+H_{\theta}+\sum_{i=1}^p\theta_iK_i
 \right),
 \\
 v_{\theta}
 &\ =\ -\sum_{i=1}^p\theta_if_i'(\alpha_i)e_i,
 \quad \alpha_i\ =\ \cos^{-1}(x^{\top}z_i),
 \quad e_i\ =\ \frac{z_i-x\cos\alpha_i}{\sin\alpha_i},
 \\
 H_{\theta}
 &\ =\ 
 e_{\theta}e_{\theta}^{\top}+\frac{\alpha_{\theta}\cos\alpha_{\theta}}{\sin\alpha_{\theta}}
 \left(I-xx^{\top}-e_{\theta}e_{\theta}^{\top}\right),
 \quad e_{\theta}\ =\ v_{\theta}/|v_{\theta}|,\quad \alpha_{\theta}\ =\ |v_{\theta}|,
 \\
 \quad K_i
 &\ =\ 
 f_i''(\alpha_i)e_ie_i^{\top}
 +\frac{f_i'(\alpha_i)\cos\alpha_i}{\sin\alpha_i}
 \left(I-xx^{\top}-e_ie_i^{\top}\right),
\end{align*}
where Euclidean coordinates in $\mbR^{n+1}$ are used.
We remark that the above formula needs
no special function,
unlike the von Mises-Fisher distribution
(\ref{eqn:vonMises-Fisher})
or the Fisher-Bingham distribution (\ref{eqn:Fisher-Bingham}).

We give examples of pairs $(f_i,z_i)$.
Recall that $W(S^n)$ is the set of all wrapping potential functions.

\begin{figure}[htbp]
  \begin{center}
   \begin{tabular}{cc}
    \includegraphics[width=6.5cm,clip]{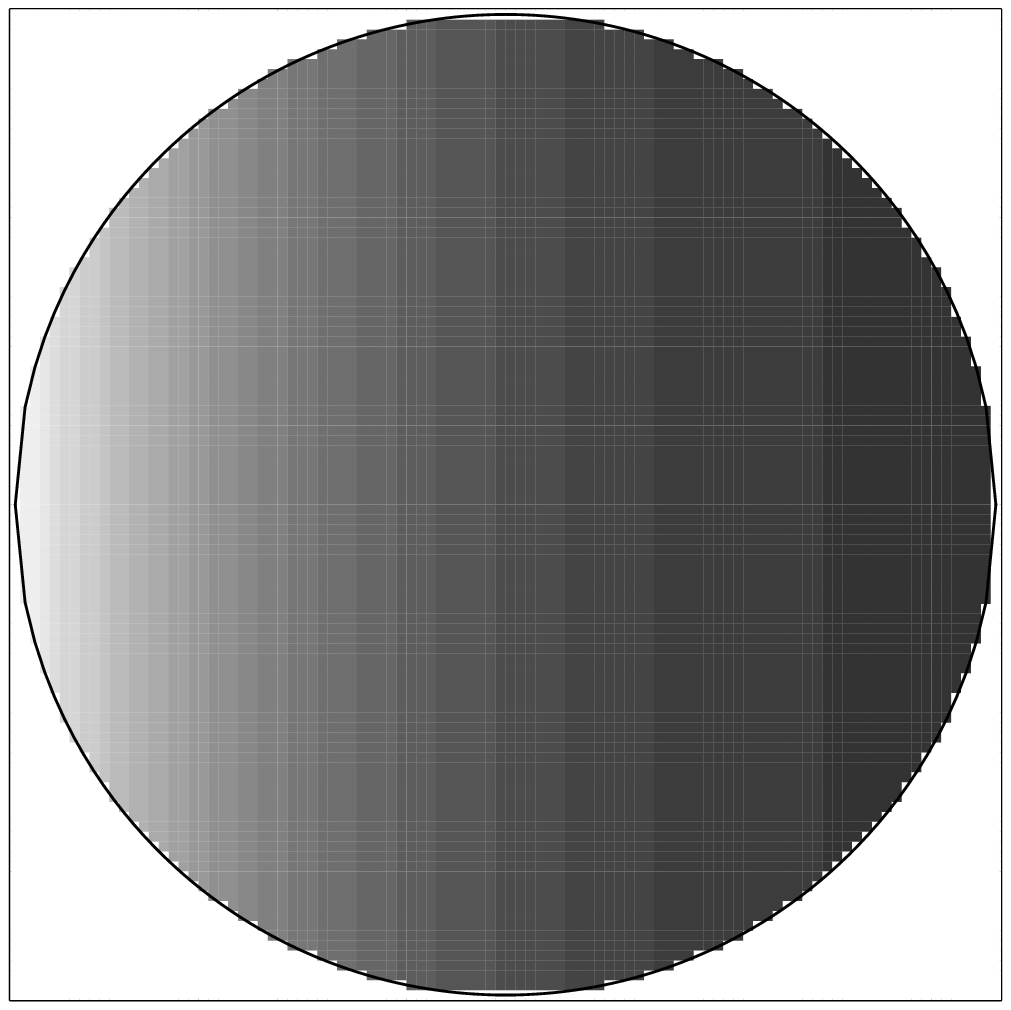} &
    \includegraphics[width=6.5cm,clip]{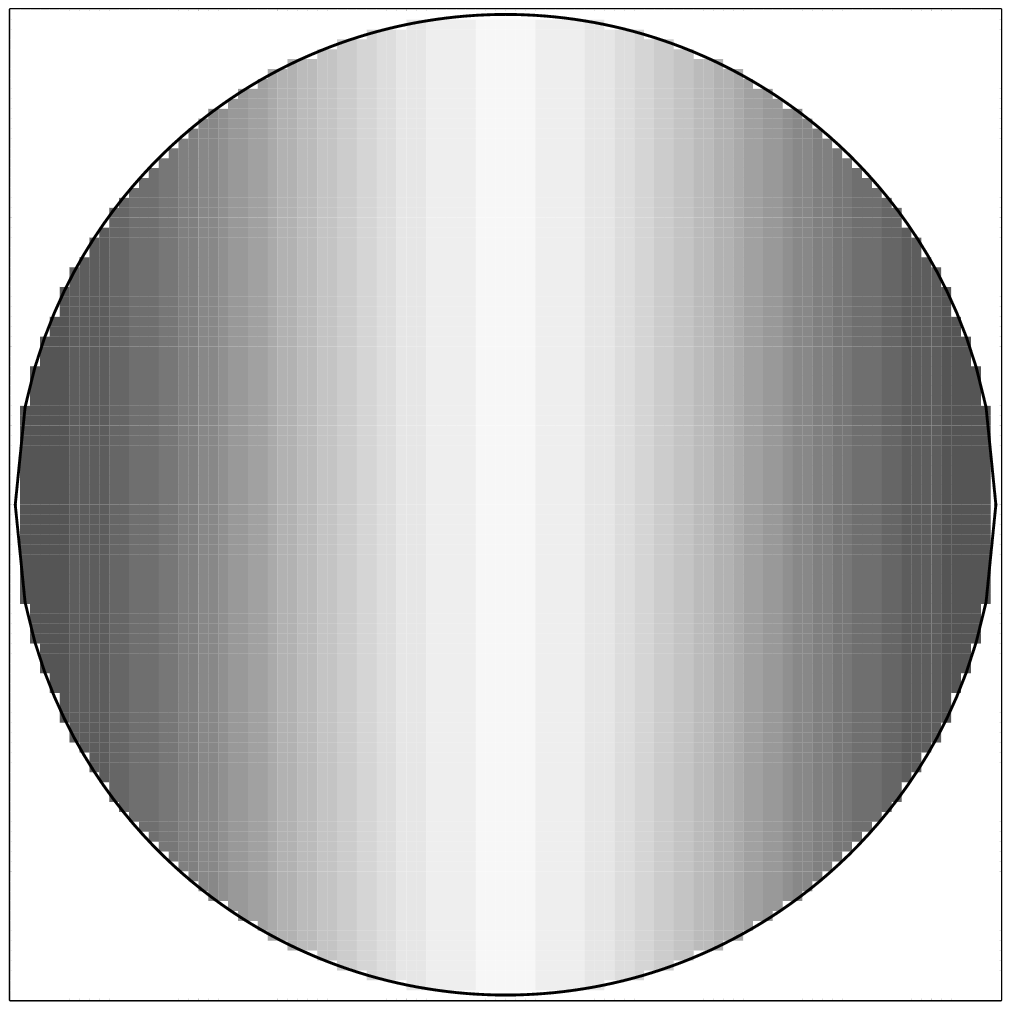}
        \vspace{-1cm}\\
    (a) $\mu=e_1$. &
    (b) $A=e_1e_1^{\top}$. \\
    \includegraphics[width=6.5cm,clip]{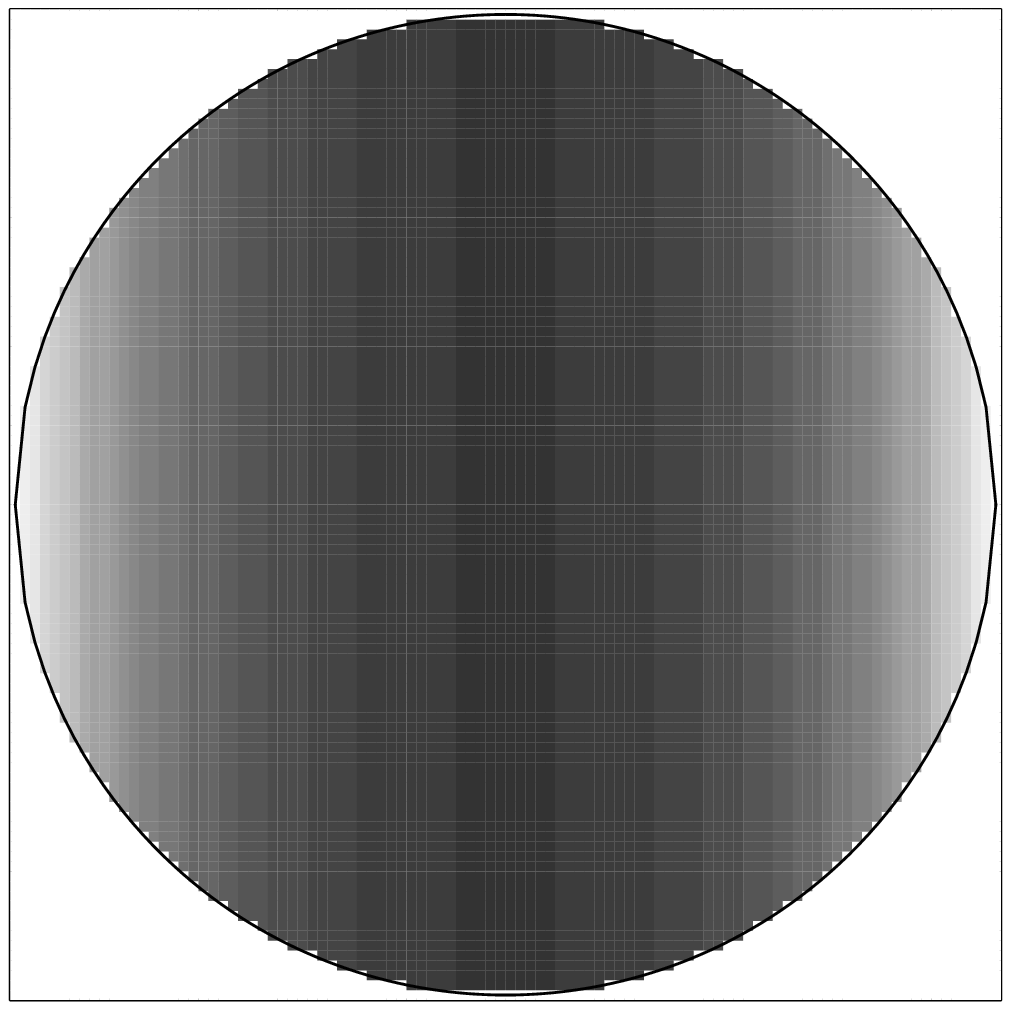} &
    \includegraphics[width=6.5cm,clip]{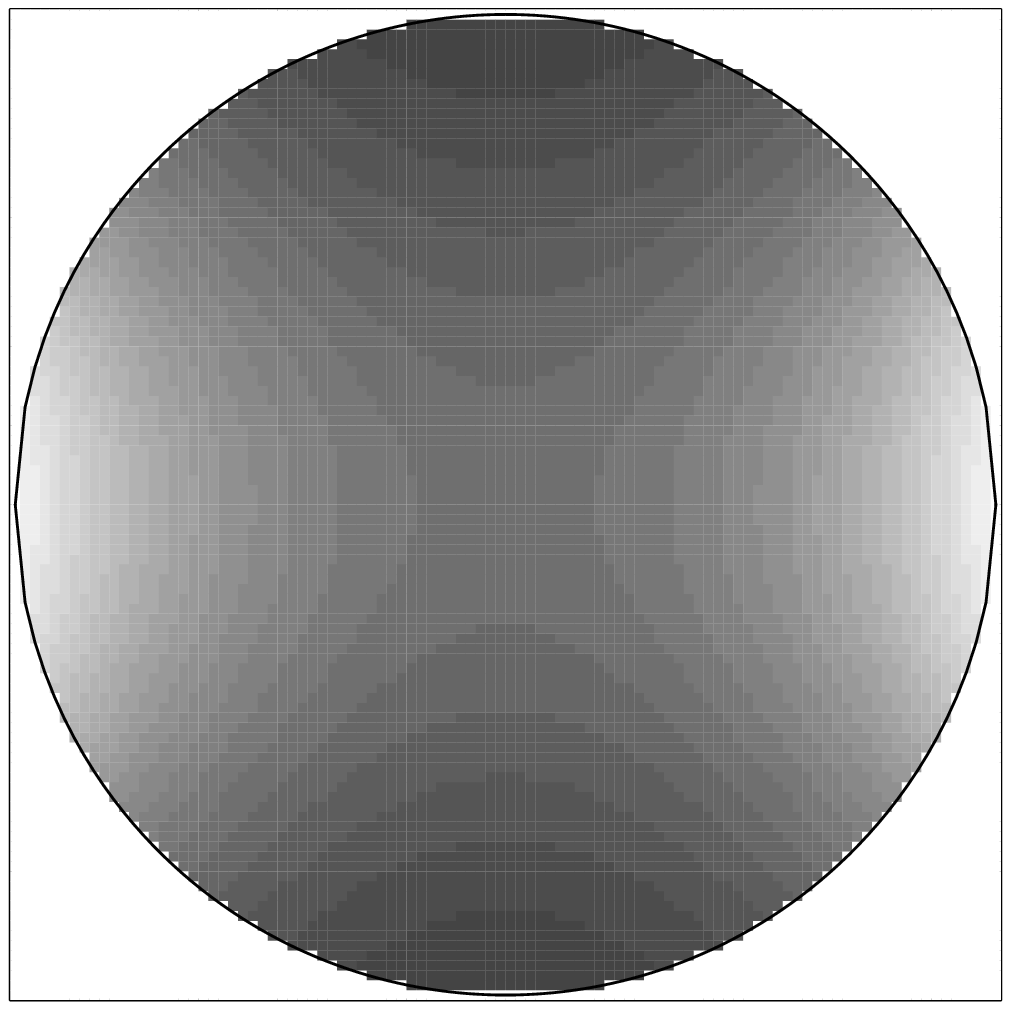}
        \vspace{-1cm}\\
    (c) $A=-e_1e_1^{\top}$. &
    (d) $A=-0.5e_1e_1^{\top}+0.5e_2e_2^{\top}$. \\
    \includegraphics[width=6.5cm,clip]{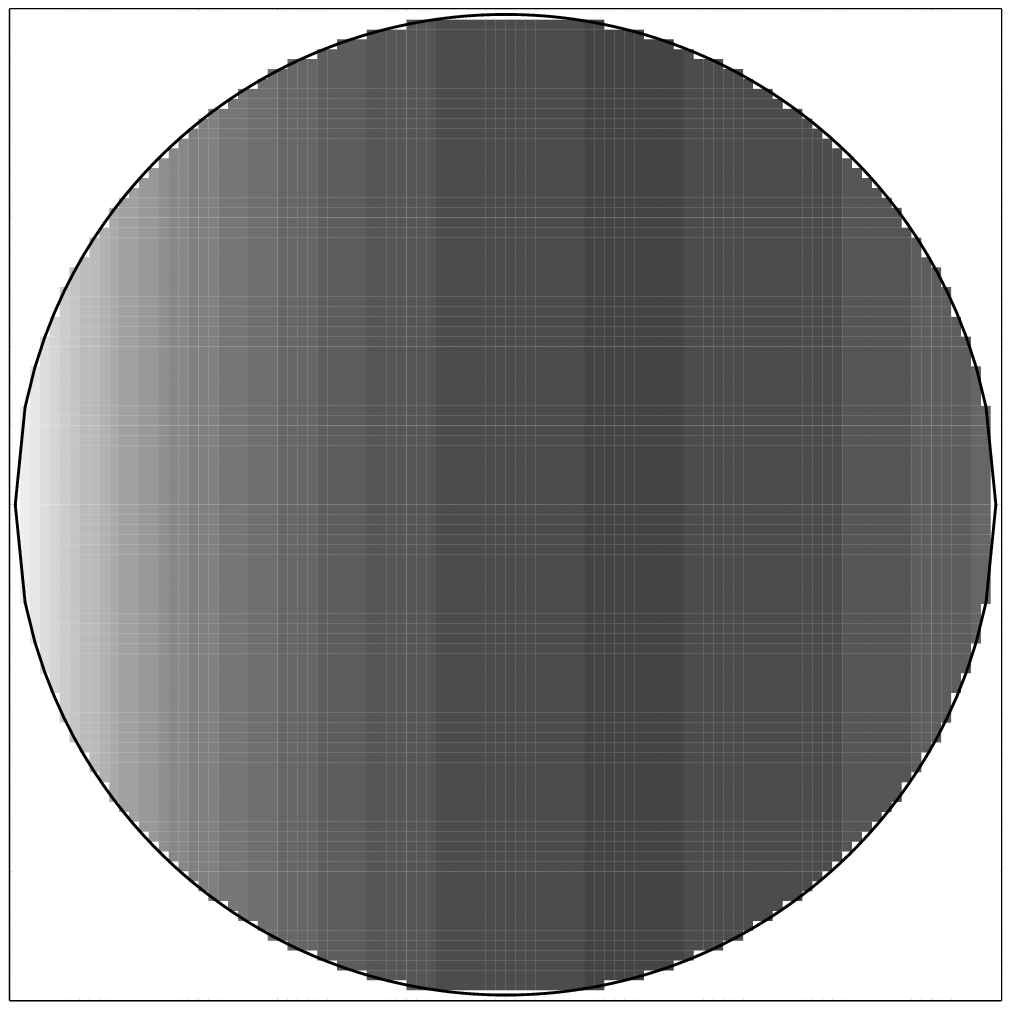} &
    \includegraphics[width=6.5cm,clip]{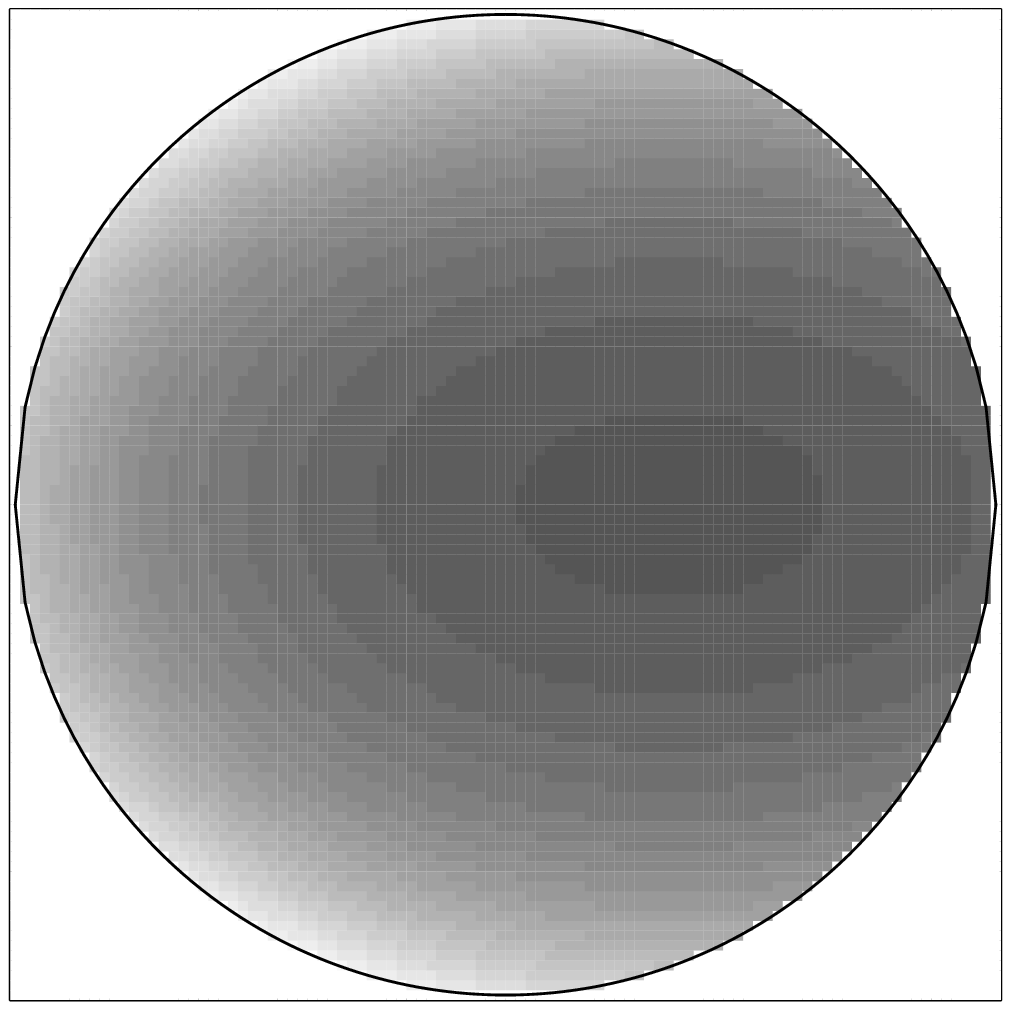}
        \vspace{-1cm}\\
    (e) $\mu=0.5e_1,A=-0.5e_1e_1^{\top}$. &
    (f) $\mu=e_1/3,A=(-e_2e_2^{\top}+e_3e_3^{\top})/3$. \\
   \end{tabular}
  \end{center}
  \caption{The quadratic-potential model.
  The white regions indicate high density.
   The figures represent (a) Concentration, (b) Negative dipole,
   (c) Positive dipole, (d) Complementary dipoles,
   (e) Unbalanced dipole and
   (f) General.
  }
   \label{fig:spherical-1}
 \end{figure}

\begin{exmp}[Linear potential]
 Let $f_i(\xi)=\cos(\xi)$ for all $i$.
 We use Euclidean coordinates in $\mbR^{n+1}$ to express $S^n$.
 Then $\phi_{\theta}(x)=\sum_{i=1}^p\theta_i\cos(d(x,z_i))=\sum_{i=1}^p \theta_ix^{\top}z_i$
 is in $W(S^n)$ as long as $\sum_{i=1}^p|\theta_i|\leq 1$.
 We deduce that a potential function $\phi_{\mu}(x):=\mu^{\top}x$
 is in $W(S^n)$ if $|\mu|\leq 1$.
 The parameter $\mu$ determines
 the direction and magnitude of concentration.
 That is, the resultant density function
 takes larger values at $x$ when $-\mu/|\mu|$ is closer to $x$
 and $|\mu|$ is larger,
 where the negative sign of $-\mu/|\mu|$ is needed because
 our model is defined by the pull-back measure.
 We call $\phi_{\mu}$ the linear potential
 and the resultant statistical model
 \textit{the linear-potential model}.
 This model is rotationally-symmetric (see Remark~\ref{rem:rotationally}).
 An example is given in Figure~\ref{fig:spherical-1} (a).
\end{exmp}

\begin{exmp}[Quadratic potential]
 Consider $f_i(\xi)=\cos(\xi)$ for $i=1,\ldots,p_1$
 and $f_i(\xi)=\cos(2\xi)/4$
 for $i=p_1+1,\ldots,p$.
 Then the potential can be written as
 \begin{align*}
  \phi_{\theta}(x)
  \ =\ \sum_{i=1}^{p_1}\theta_ix^{\top}z_i
  +\sum_{i=p_1+1}^{p}\frac{\theta_i}{4}
  \left\{2(x^{\top}z_i)^2-1\right\}.
 \end{align*}
 Let $\mu\in\mbR^{n+1}$ and $A\in\mathrm{Sym}(\mbR^{n+1})$.
 Let $|A|_1$ denote the trace norm of $A$
 defined by the sum of absolute eigenvalues of $A$.
 This is actually a norm because $|A|_1=\max_{-I\preceq B\preceq I}\mathrm{tr}[AB]$.
 Then we deduce that a potential function
 \begin{align}
  \phi_{\mu,A}(x)
  \ =\ x^{\top}\mu + \frac{1}{2}x^{\top}Ax
 \end{align}
 is in $W(S^n)$ if $(\mu,A)$ satisfies $|\mu|+|A|_1\leq 1$.
 We call the model \textit{the quadratic-potential model}.
 Various numerical examples of the quadratic potential model are
 given in Figure~\ref{fig:spherical-1}.
 Note that the representation of $A$ includes redundancy
 because $x^{\top}x=1$.
 It will be tractable if one sets $\mathrm{tr}A=0$.
 However, in general this restriction strictly reduces the size of the set.
 For example, the matrix $A=\mathrm{diag}(0.2,0,-0.8)$
 has norm $|A|_1=1$
 but the trace-adjusted one $B=\mathrm{diag}(0.4,0.2,-0.6)$
 has $|B|_1=1.2>1$.
\end{exmp}

\begin{figure}[htbp]
  \begin{center}
   \begin{tabular}{cc}
    \includegraphics[width=6.5cm,clip]{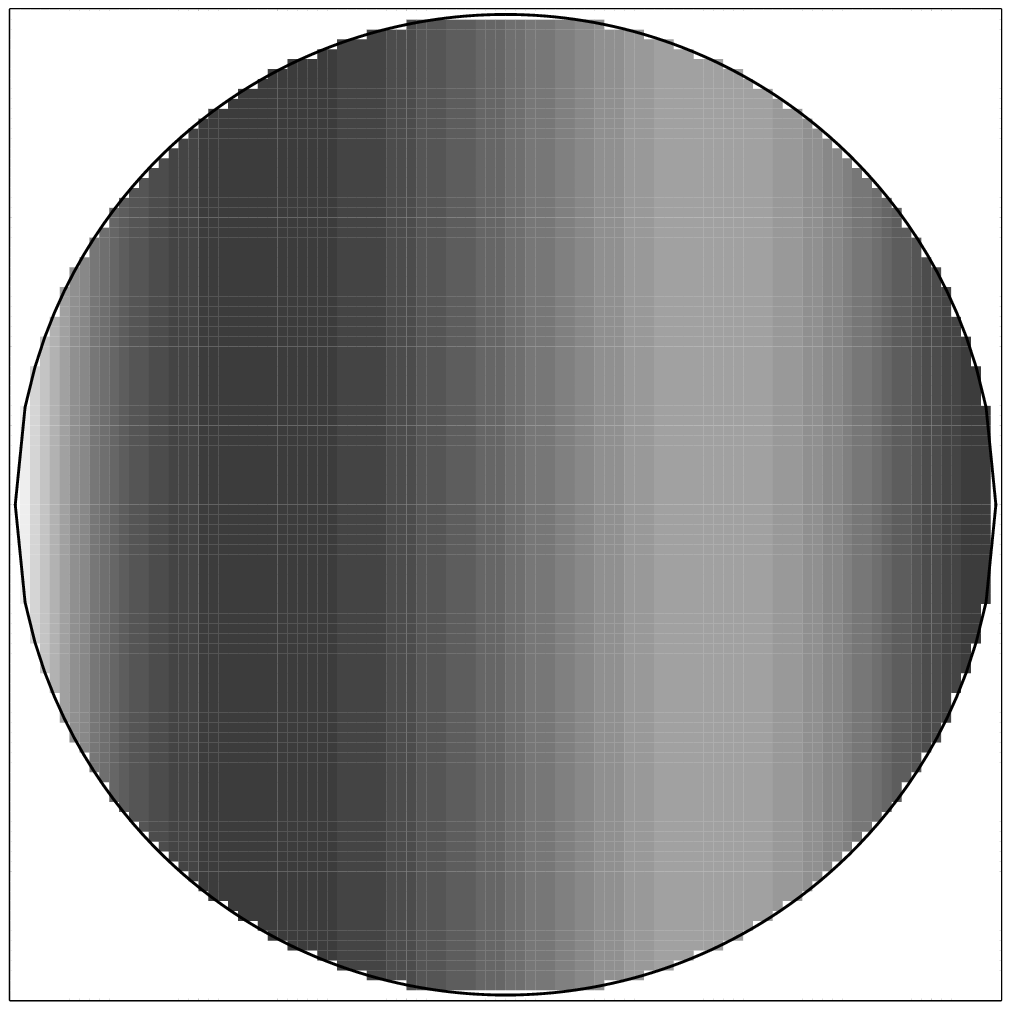} &
    \includegraphics[width=6.5cm,clip]{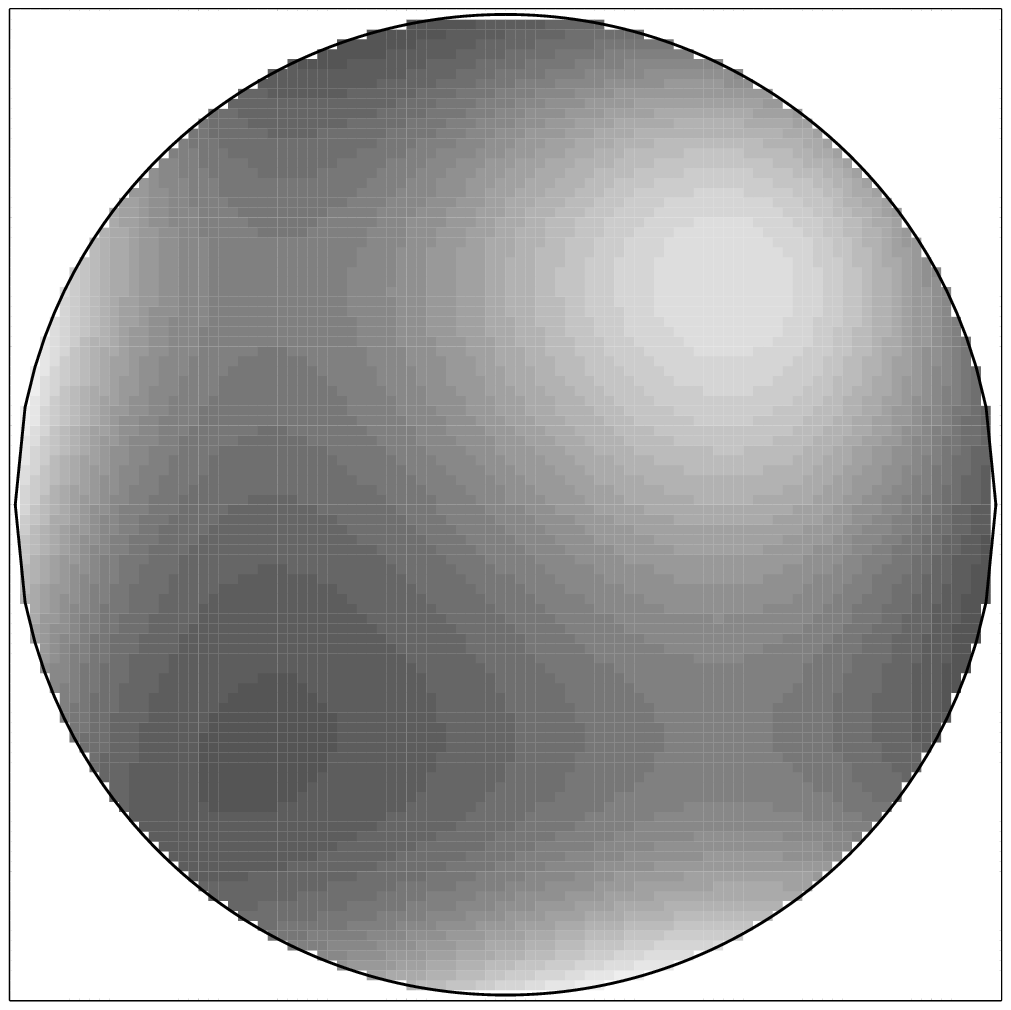}
        \vspace{-1cm}\\
    (a) $Z=(e_1),k=(3),\theta=(1)$. &
    (b) $Z=(e_1,e_2),K=(3,3),\theta=(0.5,0.5)$. \\
    \includegraphics[width=6.5cm,clip]{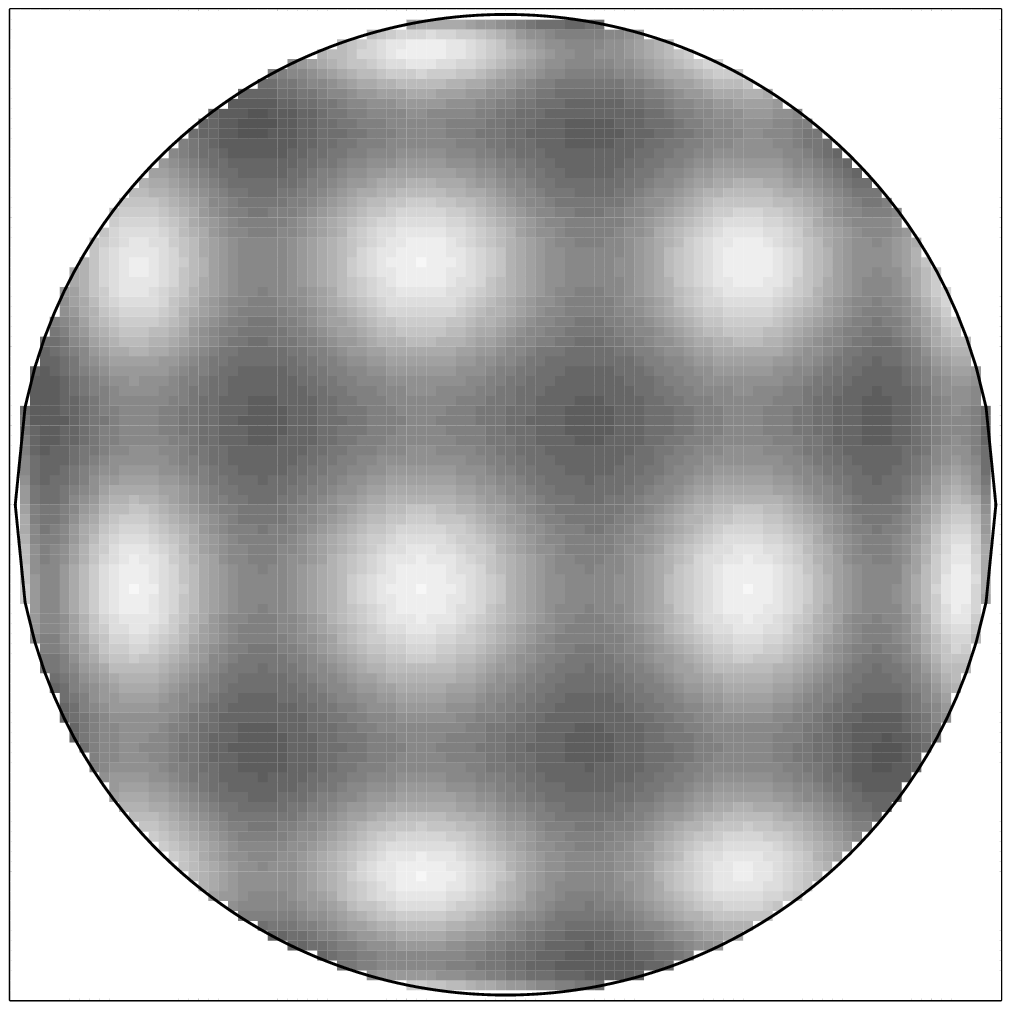} &
    \includegraphics[width=6.5cm,clip]{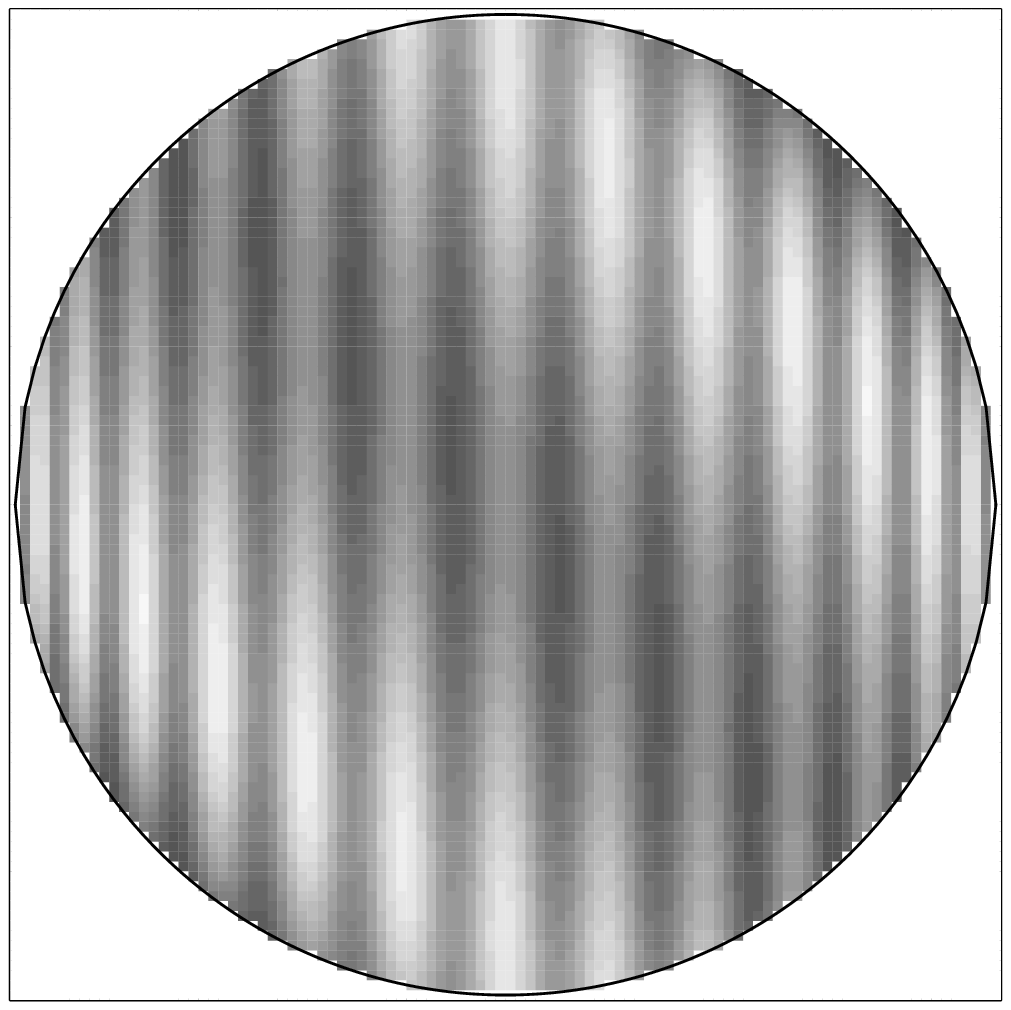}
        \vspace{-1cm}\\
    (c) $Z=(e_1,e_2),K=(9,9),\theta=(0.5,0.5)$. &
    (d) $Z=(e_1,(e_1+e_2)/\sqrt{2})$, \\
    & \quad\quad\quad $K=(30,4),\theta=(0.5,0.5)$. \\
   \end{tabular}
  \end{center}
  \caption{High-frequency spherical gradient models.
  White regions indicate high density.
  }
   \label{fig:spherical-2}
 \end{figure}

\begin{exmp}[High-frequency potential]
 As a generalization of the above examples,
 we consider $f_i(\xi)=k_i^{-2}\cos(k_i\xi)$
 for a positive integer $k_i$.
 If $Z=(z_1,\ldots,z_p)\in (S^n)^p$ and $K=(k_1,\ldots,k_p)\in\mbZ_{>0}^p$
 are given, we obtain a potential
 \begin{align}
  \phi_{\theta}(x)
  \ =\ \sum_{i=1}^p\theta_ik_i^{-2}\cos(k_id(x,z_i)).
 \end{align}
 We call this model the \textit{ high-frequency model}.
 Various numerical examples of the high-frequency model are
 given in Figure~\ref{fig:spherical-2}.
 The density function used in Figure~\ref{fig:sample}
 belongs to this class.
\end{exmp}

\subsection{An actual data set}

Here we give a brief analysis of some astronomical data.
The data consist of the locations of 188 stars of magnitude
brighter than or equal to 3.0.
The data is available from
the Bright Star Catalog (5th Revised Ed.)
distributed from the Astronomical Data Center.
We simply compare the quadratic model and the null model (uniform distribution)
by using Akaike's Information Criterion (AIC).
In general, AIC for a statistical model
is defined by the sum of $(-2)$ times the maximum log-likelihood
and $2$ times the parameter dimension.
It is recommended to select the statistical model minimizing AIC
from a set of candidates.
See \cite{akaike1974} for details of AIC.

\begin{figure}[htbp]
 \begin{center}
  \begin{tabular}{cc}
   \includegraphics[width=6.5cm,clip]{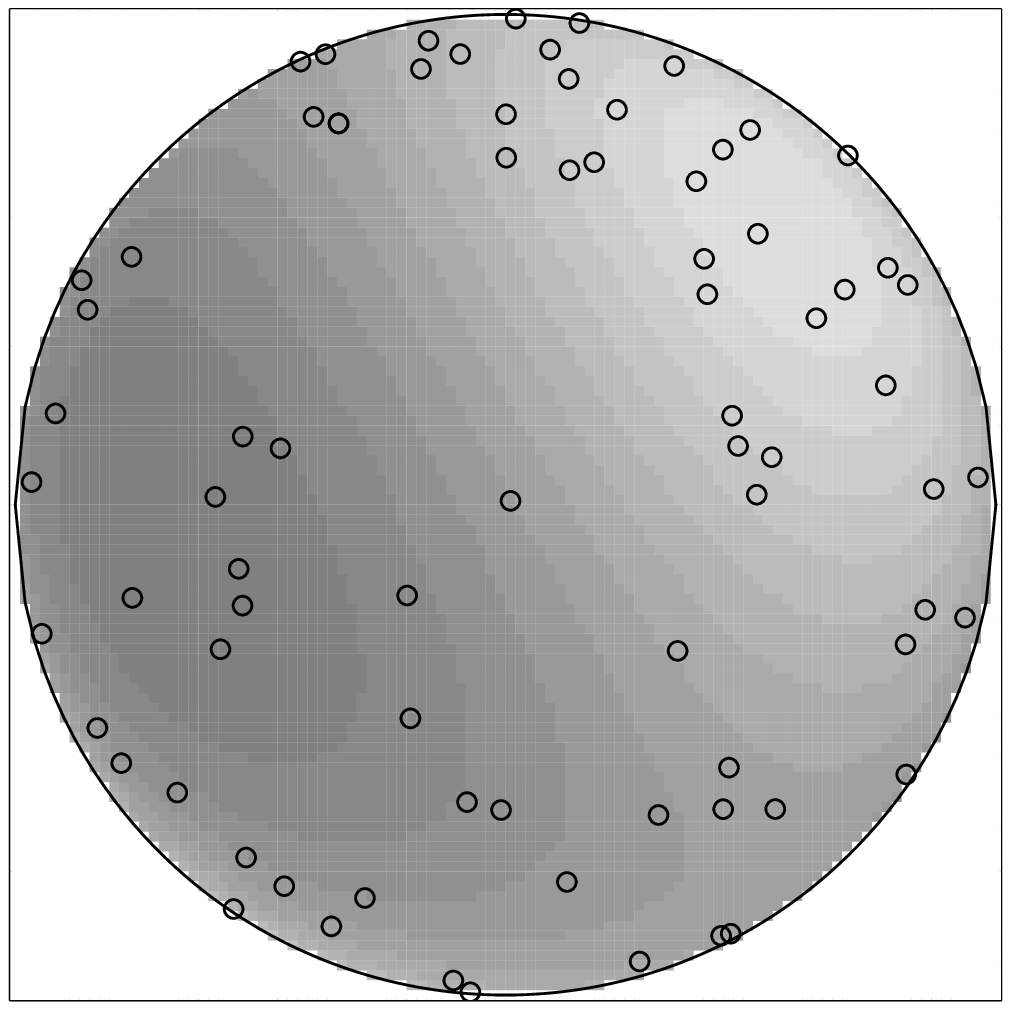} &
   \includegraphics[width=6.5cm,clip]{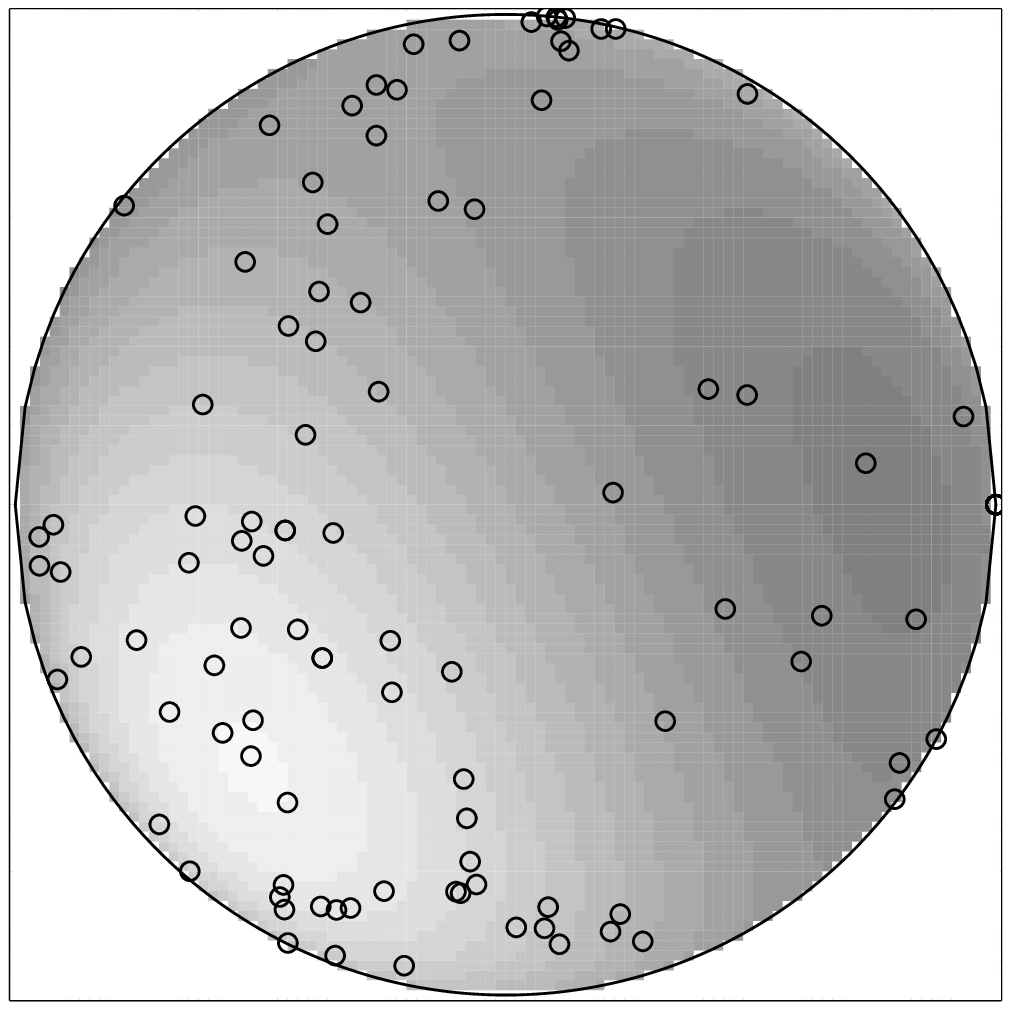}
   \vspace{-1cm}\\
   (a) Northern hemisphere. &
   (b) Southern hemisphere.
  \end{tabular}
 \end{center}
 \caption{The observed data points and the estimated density
 for the astronomic data.
 Both hemispheres are viewed from the northern side.
 White regions indicate high density.
 The points are the observed data.
 }
 \label{fig:astro}
\end{figure}

The estimated parameter for the quadratic model is
\begin{align*}
 \hat{\mu}
 \ =\ (0.010,0.017,0.091)^{\top}
 \quad \mbox{and}\quad
 \hat{A}
 \ =\ 0.173\hat{z}_1\hat{z}_1^{\top}
 -0.250\hat{z}_2\hat{z}_2^{\top},
\end{align*}
where
$\hat{z}_1
=(0.731,0.048,-0.681)^{\top}$
and $\hat{z}_2=(0.544,0.562,0.623)^{\top}$.
The maximum log-likelihood is 12.5.
Since the number of unknown parameters is 8,
AIC is $-9.0$.
On the other hand, the likelihood of the null model (uniform distribution)
is zero and AIC is also zero.
Therefore, we select the quadratic model from the two candidates.
Figure~\ref{fig:astro} shows the observed data and the estimated density.

\section{Proofs} \label{sec:proof}

\subsection{Proofs of Lemma~\ref{lem:convex}}

We use the following lemma
due to Proposition 6 of \cite{mccann2001}.
The lemma can also proved by direct calculation for the sphere.
Recall $c(x,y)=d(x,y)^2/2$.

\begin{lem}[Inverse of the exponential map] \label{lem:inverse}
 Let $x,y\in S^n$ and assume $d(x,y)<\pi$.
 Then $\nabla_xc(x,y)=-\exp_x^{-1}(y)$,
 where $\nabla_x$ denotes the gradient operator with respect to $x$.
\end{lem}

We first recall the cross-curvature non-negativity of the sphere
(\cite{kim2008}, \cite{FigalliRifford2009}).
For simplicity, the definitions below are specialized for the sphere.
For a given triplet $(x,y,z)\in (S^n)^3$ with $d(x,z)<\pi$ and $d(y,z)<\pi$,
the curve
\[
 \{\exp_z((1-t)\exp_z^{-1}(x)+t\exp_z^{-1}(y))\mid t\in[0,1]\}
\]
is called a \textit{ $c$-segment} connecting $x$ and $y$ with respect to $z$.
We denote the $c$-segment by $[x,y]_t(z)$ in this paper.
For given $x,y\in S^n$ with $d(x,y)<\pi$,
let $\sigma_s$ and $\tau_t$ be smooth curves such that
$\sigma_0=x$ and $\tau_0=y$.
We assume that either $\sigma_s=[\sigma_0,\sigma_1]_s(y)$ or $\tau_t=[\tau_0,\tau_1]_t(x)$.
Note that only one of the two curves is assumed to be a $c$-segment.
Then \textit{the cross-curvature} $\mathcal{S}$ is well defined by
\[
 \mathcal{S}(x,y)(\xi,\eta)\ =\ 
 \left.-\frac{d^2}{ds^2}\frac{d^2}{dt^2}
 c\left(\sigma_s,\tau_t\right)\right|_{s=0,t=0},
\]
where $\xi=d\sigma_s/ds|_{s=0}$ and $\eta=d\tau_t/dt|_{t=0}$.
For a given quadruplet $(x,z,y_0,y_1)\in (S^n)^4$,
\textit{the sliding mountain} is defined by a function
\begin{align}
 t\ \mapsto\ c(z,[y_0,y_1]_t(z))-c(x,[y_0,y_1]_t(z)).
 \label{eqn:sliding-mountain}
\end{align}

We use the following fact proved by \cite{kim2008} and \cite{FigalliRifford2009}.

\begin{lem}[Cross-curvature non-negativity] \label{lem:cross-curvature}
 For the sphere, the cross-curvature
 $\mathcal{S}(x,y)(\xi,\eta)$ is non-negative for any $(x,y,\xi,\eta)$
 with $d(x,y)<\pi$.
\end{lem}

Although the following lemma is essentially due to \cite{kim2008},
we derive it from Lemma~\ref{lem:cross-curvature} for completeness.

\begin{lem}[Time-convex-sliding-mountain] \label{lem:sliding-mountain}
 Let $z$ be a point in $S^n$ and
 let $y_0$ and $y_1$ be two points in $S^n$
 different from the antipodal point of $z$.
 Then for any $x\in S^n$ the sliding-mountain
 (\ref{eqn:sliding-mountain})
 is convex with respect to $t\in[0,1]$.
\end{lem}

\begin{proof}
 Denote the sliding-mountain (\ref{eqn:sliding-mountain}) by $f(t)$.
 Fix $t\in(0,1)$ and denote $y=[y_0,y_1]_t(z)$ for simplicity.
 We first assume $y$ is not the antipodal point of $x$
 and prove $d^2f(t)/dt^2\geq 0$.
 Let $\sigma_s=[z,x]_s(y)$ and $\tau_u=[y,y_1]_u(z)$.
 Note that $\sigma_s$ is a $c$-segment with respect to $\tau_0=y$.
 Then from Lemma~\ref{lem:cross-curvature}, we have
 \begin{align}
  \left.-\frac{d^2}{ds^2}\frac{d^2}{du^2}c\left(\sigma_s,\tau_u\right)
  \right|_{u=0}\ \geq\ 0
  \label{eqn:sliding-mountain-proof-1}
 \end{align}
 for each $s\in[0,1]$.
 On the other hand, by Lemma~\ref{lem:inverse},
 \begin{align*}
  \left.\frac{d}{ds}c(\sigma_s,\tau_u)\right|_{s=0}
  &\ =\ -\langle \xi,\exp_z^{-1}(\tau_u)\rangle
  \\
  &\ =\ -\langle \xi,(1-u)\exp_z^{-1}(y)+u\exp_z^{-1}(y_1)\rangle,
 \end{align*}
 where $\xi=(d\sigma_s/ds)|_{s=0}$ and $\langle\cdot,\cdot\rangle$ denotes
 the inner product on $T_zS^n$.
 We obtain
 \begin{align}
  \left.\frac{d}{ds}\frac{d^2}{du^2}c(\sigma_s,\tau_u)\right|_{s=0,u=0}
  \ =\ 0.
  \label{eqn:sliding-mountain-proof-2}
 \end{align}
 Integrating both sides of (\ref{eqn:sliding-mountain-proof-1}) with respect to $s$ twice
 and using (\ref{eqn:sliding-mountain-proof-2}), we have
 \[
  \left.\left\{\frac{d^2}{du^2}c(z,\tau_u)-
 \frac{d^2}{du^2}c(x,\tau_u)
 \right\}\right|_{u=0}
 \ \geq\ 0.
 \]
 Since $\tau_u=[y,y_1]_u(z)=[y_0,y_1]_{t+u(1-t)}(z)$, we have
 $d^2f(t)/dt^2\geq 0$.
 Next we assume that $y$ is the antipodal point of $x$.
 By assumption, $y$ is not the antipodal point of $z$.
 By direct calculation, we have
 \[
 \lim_{s\to t+0}\frac{df(s)}{ds}-\lim_{s\to t-0}\frac{df(s)}{ds}
 \ =\ 2\pi\left|\frac{d[y_0,y_1]_t(z)}{dt}\right|
 \ \geq\ 0.
 \]
 Therefore $f(t)$ is convex over $t\in[0,1]$.
\end{proof}

 Now we use Lemma~\ref{lem:sliding-mountain} to prove Lemma~\ref{lem:convex}.
 For any point $x\in S^n$, we denote the antipodal point of $x$
 by $x'$.
 Since $\phi_0$ and $\phi_1$ are $c$-convex,
 there exist functions $\phi_0^c$ and $\phi_1^c$ such that
 $\phi_i(x)=\sup_{y}\left\{-c(x,y)-\phi_i^c(y)\right\}$
 ($i=0,1$).
 Then we have
\begin{align*}
 \phi_t(x)
 &\ =\ (1-t)\phi_0(x)+t\phi_1(x)
 \\
 &\ =\ 
 \sup_{y_0}
 \left\{-(1-t)c(x,y_0)-(1-t)\phi_0^c(y_0)
 \right\}
 +\sup_{y_1}
 \left\{-tc(x,y_1)-t\phi_1^c(y_1)
 \right\}
 \\
 &\ =\ 
 \sup_{y_0\neq x'}\sup_{y_1\neq x'}
 \left\{
 -(1-t)c(x,y_0)-tc(x,y_1)
 -(1-t)\phi_0^c(y_0)-t\phi_1^c(y_1)
 \right\},
\end{align*}
 where the last equality follows from continuity of $c$ and $\phi_i^c$ ($i=0,1$).
 Now we consider a $c$-segment $[y_0,y_1]_t(z)$
 and denote it by $y_t(z)$ for simplicity.
 From Lemma~\ref{lem:sliding-mountain}, we have
\begin{align*}
 & -(1-t)c(x,y_0)-tc(x,y_1)
 \\
 &\quad\quad
 \ =\ \sup_{z\neq y_0',y_1'}
 \left\{
 -c(x,y_t(z))+c(z,y_t(z))
 -(1-t)c(z,y_0)-tc(z,y_1)
 \right\},
\end{align*}
 where the supremum of the right hand side is attained at $z=x$.
 Hence
\begin{align*}
 \phi_t(x)
 &\ =\ 
 \sup_{y_0,y_1\neq x'}\sup_{z\neq y_0',y_1'}
 \left\{
 -c(x,y_t(z))
 +c(z,y_t(z))
 -(1-t)c(z,y_0)
 -tc(z,y_1)
 \right.
 \\
 &\quad\quad\quad
 \left.
 -(1-t)\phi_0^c(y_0)-t\phi_1^c(y_1)
 \right\}
 \\
 &\ =\ \sup_{w}
 \left\{
 -c(x,w)
 -\xi(w)
 \right\},
\end{align*}
 where $\xi$ is defined by an infimum convolution
\begin{align*}
 \xi(w) &\ :=\ 
 \inf_{(y_0,y_1,z)|y_0,y_1\neq x',y_t(z)=w}
 \left\{
 -c(z,w)
 +(1-t)c(z,y_0)
 +tc(z,y_1)
 \right.
 \\
 &\quad\quad\quad
 \left.
 +(1-t)\phi_0^{c}(y_0)
 +t\phi_1^{c}(y_1)
 \right\}.
\end{align*}
 Since $\phi_t$ is written in the form
 of a $c$-transform, it is $c$-convex.
 This proves Lemma~\ref{lem:convex}.

\subsection{Proof of Theorem~\ref{thm:Jacobian}}

For each $c$-convex function $\phi$,
let $\Omega(\phi)$ be the set of points $x$
such that $|\nabla\phi(x)|<\pi$ and $\phi$ has its Hessian defined at $x$.
If $\phi$ is a wrapping potential function (Definition~\ref{defn:wrapping}),
then $S^n\setminus\Omega(\phi)$ consists only of a finite set of points.

The following lemma is essentially proved in \cite{delanoe2006}.

\begin{lem} \label{lem:less-than-pi}
 If $\phi$ is $c$-convex, then $|\nabla\phi(x)|<\pi$
 except for at most one $x\in \Omega(\phi)$.
 Furthermore, if $|\nabla\phi(x)|\geq\pi$ for some $x$,
 then $G_{\phi}(y)=G_{\phi}(x)$
 for any $y\in\Omega(\phi)$.
\end{lem}

\begin{proof}
Let $\phi$ be $c$-convex.
Assume that there exists $x\in\Omega(\phi)$
such that $|\nabla\phi(x)|\geq\pi$.
In general, any $c$-convex function on a compact Riemannian manifold
is Lipschitz continuous with Lipschitz constant less than
or equal to the diameter of the manifold (Lemma~2 of \cite{mccann2001}).
Since the diameter of the sphere $S^n$ is $\pi$,
we have $|\nabla\phi(x)|=\pi$.
Hence $G_{\phi}(x)$ is the antipodal point $x'$ of $x$.
We now prove that $G_{\phi}(y)=x'$ for all $y\in \Omega(\phi)$.
We use 2-monotonicity of the gradient map:
\begin{align*}
 d^2(x,G_{\phi}(x))+d^2(y,G_{\phi}(y))
 \ \leq\ d^2(x,G_{\phi}(y))+d^2(y,G_{\phi}(x)),
\end{align*}
where $d(x,y)$ is the distance between $x$ and $y$.
The above inequality follows from Lemma~\ref{lem:exact}.
Let $a=d(x,G_{\phi}(y))$, $b=d(y,x')$ and $c=d(y,G_{\phi}(y))$.
Then we have
$\pi^2+c^2\leq a^2+b^2$.
By the triangle inequality with respect to the triangle
$(x,y,G_{\phi}(y))$, we have $c\geq |a+b-\pi|$.
Therefore
\begin{align*}
 0
 &\ \geq\ \pi^2+c^2-a^2-b^2
 \\
 &\ \geq\ \pi^2+(a+b-\pi)^2-a^2-b^2
 \\
 &\ =\ 2(\pi-a)(\pi-b).
\end{align*}
This implies $a=\pi$ or $b=\pi$;
equivalently, $G_{\phi}(y)=x'$ or $y=x$.
Hence we have $G_{\phi}(y)=x'$ for any $y\in\Omega(\phi)$.
Then $|\nabla\phi(y)|<\pi$
for any $y\neq x$ from the definition of $G_{\phi}$.
\end{proof}

We proceed to the proof of Theorem~\ref{thm:Jacobian}.
Fix two $c$-convex functions $\phi_0$ and $\phi_1$
and let $\phi_t=(1-t)\phi_0+t\phi_1$ for $t\in[0,1]$.
Let $x\in \Omega(\phi_0)\cap\Omega(\phi_1)$.
Then it is easy to see that $x\in\Omega(\phi_t)$ for any $t\in[0,1]$.
Recall that the gradient map of $\phi_t$ is denoted by $F_t(x)=\exp_x(\nabla\phi_t(x))$.
Note that $F_t(x)$ is a $c$-segment $[F_0(x),F_1(x)]_t(x)$.
We prepare some notation to represent an explicit formula
of the Jacobian determinant of $F_t(x)$.
Let $\sigma_t(x)$ be the Jacobian determinant of the exponential map
at $\nabla\phi_t(x)$, i.e.\ 
$\sigma_t(x)=\det\{d(\exp_x(v))/dv\}|_{v=\nabla\phi_t(x)}$,
where the determinant is calculated with respect to any
orthonormal bases.
Denote the Hessian operator at $x$ by $\mathrm{Hess}_x$
and let $H_t(x)=(\mathrm{Hess}_xc(x,y))_{y=F_t(x)}$.
Then, by Cordero-Erausquin et al.~(2001),
the Jacobian determinant of $F_t(x)$ is
\begin{align}
 J_t(x)
 \ =\ 
 \sigma_t(x)
 \det
 \left(H_t(x)+\mathrm{Hess}_x\phi_t\right).
 \label{eqn:Jacobian-determinant}
\end{align}

\begin{lem} \label{lem:Hessian-concave}
 Let $x\in\Omega(\phi_0)\cap\Omega(\phi_1)$.
 The matrix-valued function $H_t(x)$
 is concave with respect to $t$:
 \[
  H_t(x)\ \succeq\ 
 (1-t)H_0(x)+t H_1(x)
 \quad \mbox{for\ any}\ t\in[0,1],
 \]
 where $A\succeq B$ means
 that $A-B$ is non-negative definite.
\end{lem}

\begin{proof}
 Since $F_t(x)$ is a $c$-segment $[F_0(x),F_1(x)]_t(x)$,
 Lemma~\ref{lem:sliding-mountain}
 implies that
\begin{align*}
 &c(w,F_t(x))-c(x,F_t(x))
 \\
 &\ \geq\ (1-t)\{c(w,F_0(x))-c(x,F_0(x))\}+t\{c(w,F_1(x))-c(x,F_1(x))\}
\end{align*}
 for all $w\in S^n$.
 By taking the Hessian with respect to $w$ at $w=x$,
 we obtain
 \[
 \mathrm{Hess}_{w}c(w,F_t(x))|_{w=x}
 \ \succeq\ (1-t)\mathrm{Hess}_wc(w,F_0(x))|_{w=x}
 +t\mathrm{Hess}_wc(w,F_1(x))|_{w=x}.
 \]
 This means $H_t(x)\succeq (1-t)H_0(x)+tH_1(x)$.
\end{proof}

\begin{lem}[Jacobian-ratio inequality] \label{lem:geometric-arithmetic}
 Let $x\in\Omega(\phi_0)\cap\Omega(\phi_1)$.
 Then the following inequality holds:
\begin{align}
 \left(
 \frac{J_t(x)}{\sigma_t(x)}
 \right)^{1/n}
 \ \geq\ 
 (1-t)
 \left(
 \frac{J_0(x)}{\sigma_0(x)}
 \right)^{1/n}
 +
 t\left(
 \frac{J_1(x)}{\sigma_1(x)}
 \right)^{1/n}.
 \label{eqn:geometric-arithmetic}
\end{align}
\end{lem}

\begin{proof}
 By the formula (\ref{eqn:Jacobian-determinant}),
 it is sufficient to prove that
 ${\det}^{1/n}(H_t+\mathrm{Hess}_x\phi_t)$
 is concave with respect to $t$.
 Indeed, by Lemma~\ref{lem:Hessian-concave} and the geometric-arithmetic inequality
 on ${\det}^{1/n}$,
 we obtain
\begin{align*}
 &{\det}^{1/n}(H_t+\mathrm{Hess}_x\phi_t)
 \\
 &\ \geq\ 
 {\det}^{1/n}\left\{
 (1-t)H_0+t H_1
 +\mathrm{Hess}_x\phi_t
 \right\}
 \\
 &\ =\ 
 {\det}^{1/n}\left\{
 (1-t)(H_0+\mathrm{Hess}_x\phi_0)
 +t(H_1+\mathrm{Hess}_x\phi_1)
\right\}
 \\
 &\ \geq\ 
 (1-t){\det}^{1/n}(H_0+\mathrm{Hess}_x\phi_0)
 +t{\det}^{1/n}(H_1+\mathrm{Hess}_x\phi_1).
\end{align*}
 Hence ${\det}^{1/n}(H_t+\mathrm{Hess}_x\phi_t)$ is concave.
\end{proof}

\begin{rem} \label{rem:Cordero}
 If $\phi_0\equiv 0$,
 the inequality (\ref{eqn:geometric-arithmetic}) is similar
 to the Jacobian inequality, due to Cordero-Erausquin et al.~(2001).
 They showed that if $\phi_0\equiv 0$,
 \begin{align}
  J_t(x)^{1/n}
 \ \geq\ (1-t)v_{1-t}(F_1(x),x)^{1/n}
 +t[v_t(x,F_1(x))]^{1/n}J_1(x)^{1/n},
  \label{eqn:Riemannian-ineq}
 \end{align}
 where $v_t(x,y)$ denotes the volume distortion
 coefficient (see Cordero-Erausquin et al.~2001 for details).
 The inequality (\ref{eqn:Riemannian-ineq}) is crucial to prove a Brunn-Minkowskii-type inequality
 on manifolds.
 However, since the inequality (\ref{eqn:Riemannian-ineq}) is
 only established for the special case $\phi_0\equiv 0$,
 it is not sufficient for our statistical application.
 Unfortunately, (\ref{eqn:Riemannian-ineq})
 is not implied from (\ref{eqn:geometric-arithmetic}).
 In fact, if $\phi_0(x)\equiv 0$, then 
 $J_0(x)=1$ and $\sigma_0(x)=1$, and
 the inequality (\ref{eqn:geometric-arithmetic}) reduces to
 \[
  J_t(x)^{1/n}
 \ \geq\ 
 (1-t)\sigma_t(x)^{1/n}
 + t\left(\frac{\sigma_t(x)}{\sigma_1(x)}\right)^{1/n}J_1(x)^{1/n}.
 \]
 This inequality is weaker than (\ref{eqn:Riemannian-ineq})
 because $v_{1-t}(F_1(x),x)>1>\sigma_t(x)$
 and $v_t(x,F_1(x))=\sigma_t(x)/\sigma_1(x)$.
 \qed
\end{rem}

\begin{lem} \label{lem:log-sigma}
 For any $x\in \Omega(\phi_0)\cap\Omega(\phi_1)$, $\log\sigma_t(x)$ is concave
 with respect to $t$.
\end{lem}

\begin{proof}
For the unit sphere $S^n$,
the Jacobian determinant of the exponential map
is given by
$(\sin|v|/|v|)^{n-1}$.
Therefore $\sigma_t(x)=(\sin|\nabla\phi_t(x)|/|\nabla\phi_t(x)|)^{n-1}$.
Since the function $[0,\pi]\ni\rho\mapsto\log(\sin\rho/\rho)$
is decreasing and concave,
and since the map $t\mapsto|\nabla\phi_t(x)|$
is convex with respect to $t$,
we deduce that the composite map
$\log\sigma_t(x)=\log(\sin|\nabla\phi_t(x)|/|\nabla\phi_t(x)|)$ is concave.
\end{proof}

 Now we prove Theorem~\ref{thm:Jacobian}.
 By Lemma~\ref{lem:geometric-arithmetic} and Lemma~\ref{lem:log-sigma},
 the functions $\log(J_t(x)/\sigma_t(x))$
 and $\log\sigma_t(x)$ are concave with respect to $t$.
 Hence $\log J_t(x)$ is also concave.

\subsection{Proof of Lemma~\ref{lem:convex-Gamma}}

Recall that $W(S^n)$ is the set of $c$-convex functions $\phi$
such that the gradient map $G_{\phi}$ is an
isomorphism on $S^n$ and $\phi$ has its Hessian defined everywhere except
at a finite set of points.

\begin{lem} \label{lem:injective}
 Let $\phi$ be a $c$-convex function and differentiable.
 Then $G_{\phi}$ is injective if and only if
 $c(x,G_{\phi}(x))+c(z,G_{\phi}(z))
 <c(x,G_{\phi}(z))+c(z,G_{\phi}(x))$
 for any $x\neq z$.
\end{lem}

\begin{proof}
 In general, by Lemma~\ref{lem:exact}, 2-monotonicity
 \[
  c(x,G_{\phi}(x))+c(z,G_{\phi}(z))
 \ \leq\ c(z,G_{\phi}(x))+c(x,G_{\phi}(z))
 \]
 holds for any $x$ and $z$, where equality holds if and only if
 $G_{\phi}(x)=G_{\phi}(z)$.
 The result follows immediately.
\end{proof}

\begin{lem} \label{lem:F_t-injective}
 Let $\phi_0$ and $\phi_1$ be members of $W(S^n)$.
 Then, for any $t\in[0,1]$,  the gradient map $F_t(x)=\exp_x(\nabla\phi_t(x))$
 is injective.
\end{lem}

\begin{proof}
 Put $h_t(x,z)=c(x,F_t(x))+c(z,F_t(z))-c(x,F_t(z))-c(z,F_t(x))$.
 By Lemma~\ref{lem:injective}, it is sufficient to show that $h_t(x,z)<0$
 for any $t\in[0,1]$ and $x\neq z$.
 By the assumption and Lemma~\ref{lem:injective},
 we have $h_0(x,z)<0$ and $h_1(x,z)<0$.
 On the other hand, by Lemma~\ref{lem:sliding-mountain},
 $h_t(x,z)\leq(1-t)h_0(x,z)+th_1(x,z)$.
 Hence we obtain $h_t(x,z)<0$.
\end{proof}

 Now we prove Lemma~\ref{lem:convex-Gamma}
 Assume that $\phi_0$ and $\phi_1$ are members of $W(S^n)$.
 From Lemma~\ref{lem:F_t-injective}, $F_t$ is injective.
 On the other hand,
 Lemma~\ref{lem:less-than-pi} implies that
 $|\nabla\phi_0(x)|<\pi$ and $|\nabla\phi_1(x)|<\pi$ 
 for all $x\in S^n$.
 Then $\nabla\phi_0(x)=\exp_x^{-1}(F_0(x))$
 and $\nabla\phi_1(x)=\exp_x^{-1}(F_1(x))$
 are continuous.
 This implies $F_t$ is continuous.
 Hence, by compactness and connectedness of $S^n$,
 $F_t$ must be an isomorphism on $S^n$.
 Twice differentiability of $\phi_t$ follows immediately
 from that of $\phi_0$ and $\phi_1$.

\subsection{Proof of Lemma~\ref{lem:rotationally-symmetric}}

 Fix $z$ and let $\phi(x)=f(d(x,z))$, for simplicity.
 We first prove $c$-convexity of $\phi$.
 It is sufficient to show that for each $x_0\in S^n$
 there exists some $y\in S^n$ such that
 $c(x_0,y)+\phi(x_0)=\inf_{x\in S^n}\{c(x,y)+\phi(x)\}$.
 Thus we investigate the point minimizing $c(x,y)+\phi(x)$
 for each fixed $y$.
 Denote the antipodal points of $y$ and $z$ by $y'$ and $z'$, respectively.
 If $x$ is different from $y'$,
 then the gradient vector of $c(x,y)+\phi(x)$
 with respect to $x$ is
 \begin{align*}
  \nabla_x\{c(x,y)+\phi(x)\}
  \ =\ \nabla_xc(x,y)
  +\nabla_xc(x,z)\frac{f^{(1)}(\sqrt{2c(x,z)})}{\sqrt{2c(x,z)}},
 \end{align*}
 where $\nabla_x$ denotes the gradient operator with respect to $x$.
 Note that the above expression makes sense for $x=z$ and $x=z'$
 because $f^{(1)}(0)=f^{(1)}(\pi)=0$ and $f\in C^2([0,\pi])$.
 By Lemma~\ref{lem:inverse},
 we know that $\nabla_xc(x,y)=-\exp_x^{-1}(y)$ and $\nabla_xc(x,z)=-\exp_x^{-1}(z)$.
 Hence the gradient vector $\nabla_x\{c(x,y)+\phi(x)\}$ vanishes only if
 $x$ lies on a great circle $C$ that passes through $y$ and $z$.
 Since the exceptional point $y'$ is also included in $C$,
 we deduce that the point minimizing $c(x,y)+\phi(x)$
 must belong to $C$.
 We fix a circular coordinate $\xi\in(-\pi,\pi]$
 representing a point on $C$
 such that $y$ corresponds to $\xi=0$.
 Let $\xi$ and $\zeta$ be the coordinates of $x$ and $z$.
 We assume $\zeta\in[0,\pi]$ without loss of generality.
 Then the function $c(x,y)+\phi(x)$ can be written as
 \begin{align*}
  h(\xi)
  \ :=\ c(x,y)+\phi(x)
  \ =\ \frac{\xi^2}{2}+f(\min\{|\xi-\zeta|,|\xi-\zeta+2\pi|\}).
 \end{align*}
 By the assumption for $f$,
 one can easily check
 that the second derivative of $h$ is $h^{(2)}(\xi)\geq 0$
 ($>0$ a.e.)
 as long as
 $\xi\neq \pi$.
 Furthermore,
 we obtain $h^{(1)}(\pi-0)>h^{(1)}(-\pi+0)$.
 Thus $\xi=\pi$ is not a point minimizing $h$.
 Furthermore, the point minimizing $h$ is unique because
 $h$ is strictly convex over $(-\pi,\pi)$.
 We denote the minimizer by $\xi_0\in (-\pi,\pi]$
 and the corresponding point in $S^n$ by $x_0$.
 If $y$ revolves along a great circle $C$
 passing through $z$,
 then $x_0$ must continuously revolve along $C$.
 Since $y$ can belong any great circle passing through $z$,
 we deduce that for each point $x_0$
 there exists some $y\in C$ such that
 the function $c(x,y)+\phi(x)$ of $x\in S^n$ is minimized at $x_0$.
 This proves $c$-convexity of $\phi$.

 Next we prove the gradient map $G_{\phi}(x)$ is well defined and an isomorphism.
 Since $\phi$ is differentiable everywhere,
 $G_{\phi}$ is well defined.
 Let $x=\exp_z(t e)$, where $t\in[0,\pi]$ and $e\in T_zS^n$ with $|e|=1$.
 Then the gradient map is explicitly given by
  $G_{\phi}(x)=\exp_z\left((t+f^{(1)}(t))e\right)$.
 If $t$ moves from $0$ to $\pi$,
 then $t+f^{(1)}(t)$ moves from $0$ to $\pi$ monotonically
 because $1+f^{(2)}(t)>0$ for almost all $t\in[0,\pi]$.
 Hence $G_{\phi}:S^n\to S^n$ is an isomorphism.

 Lastly, $\phi$ is clearly twice differentiable
 whenever $x\neq z$ and $x\neq z'$.
 This completes the proof.

\section{Discussion} \label{sec:discussion}

We briefly discuss the Jacobian inequality
for general manifolds.

In the proof of Theorem~\ref{thm:Jacobian},
we have used
the closed property of cost-convex functions (Lemma~\ref{lem:convex}),
the Jacobian-ratio inequality (Lemma~\ref{lem:geometric-arithmetic})
and log-concavity of the Jacobian of the exponential map
(Lemma~\ref{lem:log-sigma}).
For any non-negatively cross-curved
(or time-convex-sliding-mountain) manifold defined in \cite{kim2008},
the former two lemmas are obtained in the same manner.
However, Lemma~\ref{lem:log-sigma} does not automatically follow from
the non-negative cross-curvature condition.

The author does not know
if any Riemannian manifold with non-negative cross-curvature
satisfies the Jacobian inequality.
At least, any product space of $S^n$ and $\mbR^n$
satisfies the Jacobian inequality
because the non-negative cross-curvature condition
is preserved for products of manifolds (\cite{kim2008})
and the Jacobian determinant of the exponential map is also
factorized into the Jacobian determinant on each space.
This fact may enable us to describe dependency structures
of multivariate directional data in statistics.
We leave such an extension for future research.

\section*{Acknowledgements}

The author is grateful to Alessio Figalli and Robert J.\ McCann for
their helpful comments on the first version of the paper.
This study was partially supported by the Global Center of
Excellence ``The research and training center for new development in
mathematics''
and by the Ministry of Education, Science, Sports and Culture, Grant-in-Aid
for Young Scientists (B), No.~19700258.



\newpage

\end{document}